\documentclass[11pt]{amsart}
\usepackage{amssymb}
\usepackage{amsfonts}
\usepackage{amsmath}

\newtheorem{theorem}{Theorem}[section]

\newtheorem{lemma}[theorem]{Lemma}
\theoremstyle{definition}

\newtheorem{acknowledgement}{Acknowledgement}

\theoremstyle{remark}
\newtheorem{remark}[theorem]{Remark}
\numberwithin{equation}{section}

\begin{document}
\title{Nodal solutions for singular semilinear elliptic systems}
\author[A. Moussaoui]{Abdelkrim Moussaoui}
\address{A. Moussaoui\\
Applied Mathematics Laboratory (LMA), Faculty of Exact Sciences\\
and Biology departement, Faculty of Natural \& Life Sciences\\
A. Mira Bejaia University, Targa Ouzemour, 06000 Bejaia, Algeria}
\email{abdelkrim.moussaoui@univ-bejaia.dz}

\begin{abstract}
In this paper, we prove existence of nodal solutions for singular semilinear
elliptic systems without variational structure where its both components are
of sign changing. Our approach is based on sub-supersolutions method
combined with perturbation arguments involving singular terms.
\end{abstract}

\maketitle

\section{Introduction}

\label{S1}

Let $\Omega $ be a bounded domain in $\mathbb{R}^{N}$ ($N\geq 2)$ having a
smooth boundary $\partial \Omega $ and a positive measure such that $%
meas(\Omega )>1$. Consider the following system of semilinear elliptic
equations%
\begin{equation*}
(\mathrm{P})\qquad \left\{ 
\begin{array}{ll}
-\Delta u+h_{\lambda ,\phi _{1}}(u)=a_{1}(x)\frac{f_{1}(v)}{|u|^{\alpha _{1}}%
} & \text{in }\Omega \\ 
-\Delta v+h_{\lambda ,\phi _{1}}(v)=a_{2}(x)\frac{f_{2}(u)}{|v|^{\alpha _{2}}%
} & \text{in }\Omega \\ 
u,v=0 & \text{on }\partial \Omega ,%
\end{array}%
\right.
\end{equation*}%
where $\Delta $ stands for the Laplace differential operator and $h_{\lambda
,\phi _{1}}$ is a linear function defined by%
\begin{equation}
h_{\lambda ,\phi _{1}}(s):=\lambda (s+\phi _{1}),\text{ for }s\in 
\mathbb{R}
,\text{ for }\lambda >0,  \label{24}
\end{equation}%
where $\phi _{1}$\ denotes the positive eigenfunction corresponding to the
principal eigenvalue $\lambda _{1}$. In the reaction terms (the right hand
side) of problem $(\mathrm{P}),$ the function $a_{i}\in L^{\infty }(\Omega )$
satisfies

\begin{description}
\item[$\mathrm{H(a)}$] There exists a constant $1<\rho _{i}<meas(\Omega )$
such that 
\begin{equation*}
\left\{ 
\begin{array}{l}
a_{i}(x)>0\ \text{for a.a. }x\in \Omega _{\rho _{i}} \\ 
a_{i}(x)\leq 0\ \text{for a.e. }x\in \Omega \backslash \overline{\Omega }%
_{\rho _{i}},%
\end{array}%
\right.
\end{equation*}%
where%
\begin{equation*}
\Omega _{\rho _{i}}=\left\{ x\in \Omega :d\left( x,\partial \Omega \right)
<\rho _{i}\right\} ,\text{ for }i=1,2,
\end{equation*}
\end{description}

while the nonlinear term $f_{i}$ is a continuous function satisfying the
growth condition

\begin{description}
\item[$\mathrm{H}(f)$] There is constants $m_{i},M_{i}>0$ and $\beta _{i}\in
(0,1)$ such that%
\begin{equation*}
m_{i}\leq f_{i}(s)\leq M_{i}(1+|s|^{\beta _{i}}),\text{ \ for all }s\in 
\mathbb{R}
.
\end{equation*}
\end{description}

We consider the system $(\mathrm{P})$ in a singular case assuming that%
\begin{equation}
0<\alpha _{1},\alpha _{2}<1.  \label{exp}
\end{equation}

\mathstrut

Our main goal is to provide a nodal solution $(u,v)$ of singular
nonvariational elliptic system $(\mathrm{P})$. This means that both
components $u$ and $v$ are of sign changing. According to our knowledge,
this topic is a novelty. The virtually non-existent works in the literature
devoted to this subject is partly due to the fact that the existence of a
nodal solution for systems is more delicate than in the case of equations.
Precisely, \cite{MMP} is the only paper that has addressed this issue for
nonvariational systems through topological degree argument where a different
concept of nodal solutions is introduced. Namely, it is shown the existence
of a solution where its both components$\ $are nontrivial and are not of the
same constant sign. Solutions of this type have also been studied in \cite%
{M, MZ} for a class of quasilinear systems with variational structure by
combining variational methods with suitable truncation. However, the nodal
solution in \cite{M} is defined in more subtle way considering that either
its components are of the same constant sign or at least one of them is sign
changing. Thence, even for variational systems, the question regarding nodal
solutions where both components are of sign changing remains open. Here, it
should be noted that system $(\mathrm{P})$ under assumptions above is not in
variational form, so the variational methods are not applicable.

\mathstrut

Another main technical difficulty in the present paper consists in the
presence of singularities in system $(\mathrm{P})$ near the origin that
occur under assumption (\ref{exp}). These singularities make difficult any
study attendant to nodal solutions for $(\mathrm{P})$ due to their sign
change forcing them to cross zero. It represents a serious difficulty to
overcome and, as far as we know, is never handled in the literature even for
singular problems in the scalar case. For more inquiries on the study of
constant sign solutions for singular systems we refer to \cite{MM1, MM2,
MM3, DM1, DM2, KM} and the references therein.

To handle our problem, we show that the sets where $u$ and $v$ vanish are of
zero measure. This is an essential point enabling nodal solutions
investigation. Thereby, by a solution of problem $(\mathrm{P})$ we mean a
couple $(u,v)\in H_{0}^{1}\left( \Omega \right) \times H_{0}^{1}\left(
\Omega \right) $ such that the set where $u$ (resp. $v$) vanishes is
negligible and%
\begin{equation*}
\left\{ 
\begin{array}{c}
\int_{\Omega }(\nabla u\nabla \varphi +h_{\lambda ,\phi _{1}}(u)\varphi )\
dx=\int_{\Omega }a_{1}(x)\frac{f_{1}(v)}{|u|^{\alpha _{1}}}\varphi \ dx \\ 
\int_{\Omega }(\nabla v\nabla \psi +h_{\lambda ,\phi _{1}}(v)\psi )\
dx=\int_{\Omega }a_{2}(x)\frac{f_{2}(u)}{|v|^{\alpha _{2}}}\psi \ dx,%
\end{array}%
\right.
\end{equation*}%
for all $(\varphi ,\psi )\in H_{0}^{1}\left( \Omega \right) \times
H_{0}^{1}\left( \Omega \right) .$

\mathstrut

Our approach is chiefly based on sub-supersolution method. It is applied to
a disturbed system $(\mathrm{P}_{\varepsilon })$ depending on parameter $%
\varepsilon >0$ whose study is relevant for problem $(\mathrm{P})$. The
obtained solution $(u_{\varepsilon },v_{\varepsilon })$ of $(\mathrm{P}%
_{\varepsilon })$ is located in the rectangle formed by sub-supersolutions.
A significant feature of our result lies in the construction of the sub- and
supersolution pair for $(\mathrm{P}_{\varepsilon })$. At this point, the
choice of suitable functions with an adjustment of adequate constants is
crucial. Namely, exploiting spectral properties of the Laplacian operator,
the supersolution $(\overline{u},\overline{v})$, constructed explicitly, is
sign-changing and independent of $\varepsilon >0,$ while the subsolution $(%
\underline{u}_{\varepsilon },\underline{v}_{\varepsilon })$ which, besides
the dependence of $\varepsilon ,$ does not have an explicit form, admits a
limit as $\varepsilon \rightarrow 0$ the couple $(\underline{u},\underline{v}%
)$ where the component $\underline{u}$ (resp. $\underline{v}$) is negative
in $\Omega _{\rho _{1}}$ (resp. $\Omega _{\rho _{2}})$ and nonnegative (not
necessarly positive) in $\Omega \backslash \overline{\Omega }_{\rho _{1}}$
(resp. $\Omega \backslash \overline{\Omega }_{\rho _{2}}$). Actually, it is
worth noting that $(\underline{u}_{\varepsilon },\underline{v}_{\varepsilon
})$ is a solution of an auxiliary problem $(\mathrm{\tilde{P}}_{\varepsilon
})$ related to $(\mathrm{P}_{\varepsilon })$. Then, the general theory of
sub-supersolutions for systems of quasilinear equations (see \cite{CLM})
implies the existence of a solution $(u_{\varepsilon },v_{\varepsilon })$ of
problem $(\mathrm{P}_{\varepsilon })$ with the sets where $u_{\varepsilon }$
and $v_{\varepsilon }$ vanish are negligible. In particular, this establish
that $u_{\varepsilon }$ and $v_{\varepsilon }$ cannot be identically zero in 
$\Omega _{\rho _{1}}$ and $\Omega _{\rho _{2}}$, respectively. Then, the
solution $(u,v)$ of $(\mathrm{P}),$ lying in $[\underline{u},\overline{u}%
]\times \lbrack \underline{v},\overline{v}]$ with the property that $u\neq 0$
in $\Omega _{\rho _{1}}$ and $v\neq 0$ in $\Omega _{\rho _{2}}$, is derived
by passing to the limit as $\varepsilon \rightarrow 0.$ The argument is
based on a priori estimates, dominated convergence Theorem as well as $S_{+}$%
-property of the negative Laplacian. Hence, $(u,v)$ turns out a nodal
solution of $(\mathrm{P})$ with its both components $u$ and $v$ are of sign
changing.

\mathstrut

The rest of this article is organized as follows. Section \ref{S2} contains
the proof of the existence of solutions for regularized system $(\mathrm{P}%
_{\varepsilon })$ as well as the construction of its sign changing
sub-supersolutions. Section \ref{S3} presents the proof of the existence of
nodal solutions of system $(\mathrm{P})$.

\section{The regularized system}

\label{S2}

For $\varepsilon >0$, let consider the auxiliary system 
\begin{equation*}
(\mathrm{P}_{\varepsilon })\qquad \left\{ 
\begin{array}{ll}
-\Delta u+h_{\lambda ,\phi _{1}}(u)=a_{1}(x)\frac{f_{1}(v)}{(|u|+\varepsilon
)^{\alpha _{1}}} & \text{in }\Omega \\ 
-\Delta v+h_{\lambda ,\phi _{1}}(v)=a_{2}(x)\frac{f_{2}(u)}{(|v|+\varepsilon
)^{\alpha _{2}}} & \text{in }\Omega \\ 
u,v=0 & \text{on }\partial \Omega .%
\end{array}%
\right.
\end{equation*}%
Employing sub-supersolution method we shall prove that problem $(\mathrm{P}%
_{\varepsilon })$ admits a nontrivial solution.

We recall that a sub-supersolution for $(\mathrm{P}_{\varepsilon })$ is any
pairs $(\underline{u},\underline{v})\in (H_{0}^{1}\left( \Omega \right) \cap
L^{\infty }(\Omega ))^{2}$ and $(\overline{u},\overline{v})\in (H^{1}\left(
\Omega \right) \cap L^{\infty }(\Omega ))^{2}$ for which there hold $(%
\overline{u},\overline{v})\geq (\underline{u},\underline{v})$ a.e. in $%
\overline{\Omega }$, 
\begin{equation*}
\left\{ 
\begin{array}{l}
\int_{\Omega }(\nabla \underline{u}\nabla \varphi +h_{\lambda ,\phi _{1}}(%
\underline{u})\varphi )\ \mathrm{d}x-\int_{\Omega }a_{1}(x)\frac{f_{1}(v)}{(|%
\underline{u}|+\varepsilon )^{\alpha _{1}}}\varphi \ \mathrm{d}x\leq 0 \\ 
\int_{\Omega }(\nabla \underline{v}\nabla \psi +h_{\lambda ,\phi _{1}}(%
\underline{v})\psi )\ \mathrm{d}x-\int_{\Omega }a_{2}(x)\frac{f_{2}(u)}{(|%
\underline{v}|+\varepsilon )^{\alpha _{2}}}\psi \ dx\leq 0,%
\end{array}%
\right.
\end{equation*}%
\begin{equation*}
\left\{ 
\begin{array}{l}
\int_{\Omega }(\nabla \overline{u}\nabla \varphi +h_{\lambda ,\phi _{1}}(%
\overline{u})\varphi )\ \mathrm{d}x-\int_{\Omega }a_{1}(x)\frac{f_{1}(v)}{(|%
\overline{u}|+\varepsilon )^{\alpha _{1}}}\varphi \ \mathrm{d}x\geq 0 \\ 
\int_{\Omega }(\nabla \overline{v}\nabla \psi +h_{\lambda ,\phi _{1}}(%
\overline{v})\psi )\ \mathrm{d}x-\int_{\Omega }a_{2}(x)\frac{f_{2}(u)}{(|%
\overline{v}|+\varepsilon )^{\alpha _{2}}}\psi \ \mathrm{d}x\geq 0,%
\end{array}%
\right.
\end{equation*}%
for all $\varphi ,\psi \in H_{0}^{1}\left( \Omega \right) $ with $\varphi
,\psi \geq 0$ a.e. in $\Omega $ and for all $u,v\in H_{0}^{1}\left( \Omega
\right) $ satisfying $\underline{u}\leq u\leq \overline{u}$ and $\underline{v%
}\leq v\leq \overline{v}$ a.e. in $\Omega $.

\subsection{A constant sign sub-supersolution pair}

In what follows $\phi _{1}$\ denotes the positive eigenfunction
corresponding to the principal eigenvalue $\lambda _{1}$, that is,%
\begin{equation*}
-\Delta \phi _{1}=\lambda _{1}\phi _{1}\text{ in }\Omega ,\text{ }\phi _{1}=0%
\text{ on }\partial \Omega .
\end{equation*}%
which is well known to verify 
\begin{equation}
\begin{array}{c}
l^{-1}d(x)\leq \phi _{1}(x)\leq ld(x)\text{ for all }x\in \Omega ,%
\end{array}
\label{67}
\end{equation}%
\begin{equation}
\begin{array}{l}
|\nabla \phi _{1}|\geq \eta \text{ \ as }d(x)\rightarrow 0,%
\end{array}
\label{19}
\end{equation}%
with constants $l>1$ and $\eta >0$, where $d(x)$ denotes the distance from a
point $x\in \overline{\Omega }$ to the boundary $\partial \Omega $ and $%
\overline{\Omega }=\Omega \cup \partial \Omega $ is the closure of $\Omega
\subset 
\mathbb{R}
^{N}$.

Let $\tilde{\Omega}$ be a bounded domain in $%
\mathbb{R}
^{N}$ with a smooth boundary $\partial \tilde{\Omega}$ such that $\overline{%
\Omega }\subset \tilde{\Omega}.$ Denote $\tilde{d}(x):=d(x,\partial \tilde{%
\Omega})$. By the definition of $\tilde{\Omega}$ there exists a constant $%
\mu >0$ small enough such that 
\begin{equation}
\tilde{d}(x)>\mu \text{ in }\overline{\Omega }.  \label{18}
\end{equation}%
Let $\tilde{e}\in C^{1}(\overline{\tilde{\Omega}})$ be the unique solution
of the Dirichlet problem%
\begin{equation}
\left\{ 
\begin{array}{ll}
-\Delta \tilde{e}=1 & \text{in }\tilde{\Omega} \\ 
\tilde{e}=0 & \text{on }\partial \tilde{\Omega},%
\end{array}%
\right.  \label{3*}
\end{equation}%
which is known to satisfy the estimate 
\begin{equation}
c^{-1}\tilde{d}(x)\leq \tilde{e}(x)\leq c\tilde{d}(x)\ \text{\ in\ }\tilde{%
\Omega},  \label{3}
\end{equation}%
for certain constant $c>1$ (see \cite[proof of Lemma 3.1]{DM}).

\begin{theorem}
\label{T2} Under assumptions $\mathrm{H}(f),$ $\mathrm{H}(a)$ and (\ref{exp}%
), system $(\mathrm{P}_{\varepsilon })$ possesses a solution $%
(u_{\varepsilon },v_{\varepsilon })\in (H_{0}^{1}(\Omega )\cap L^{\infty
}(\Omega ))^{2}$ within $[-C\tilde{e},C\tilde{e}]\times \lbrack -C\tilde{e},C%
\tilde{e}],$ with a constant $C>1$ large, for all $\varepsilon \in (0,1)$
and all $\lambda \geq 0$.
\end{theorem}

\begin{proof}
Using (\ref{18})-(\ref{3}) furnishes%
\begin{equation}
\begin{array}{l}
-\Delta (C\tilde{e})+h_{\lambda ,\phi _{1}}(C\tilde{e})=C+\lambda (C\tilde{e}%
+\phi _{1}) \\ 
\geq C(1+\lambda \tilde{e})\geq C(1+\lambda \frac{\mu }{c})\text{ \ in }%
\Omega%
\end{array}
\label{61}
\end{equation}%
and%
\begin{equation}
\begin{array}{l}
-\Delta (-C\tilde{e})+h_{\lambda ,\phi _{1}}(-C\tilde{e})=-C(1+\lambda 
\tilde{e})+\lambda \phi _{1} \\ 
\leq -C(1+\lambda \frac{\mu }{c})+\lambda \left\Vert \phi _{1}\right\Vert
_{\infty }\text{ in }\Omega .%
\end{array}
\label{62}
\end{equation}%
By $\mathrm{H}(f),$ $\mathrm{H}(a),$ (\ref{exp}) and (\ref{18}), it hold%
\begin{equation}
\begin{array}{l}
a_{1}(x)\frac{f_{1}(u)}{(|C\tilde{e}|+\varepsilon )^{\alpha _{1}}}\leq
M_{1}\left\Vert a_{1}\right\Vert _{\infty }\frac{1+|v|^{\beta _{1}}}{(|C%
\tilde{e}|+\varepsilon )^{\alpha _{1}}}\leq M_{1}\left\Vert a_{1}\right\Vert
_{\infty }\frac{1+|C\tilde{e}|^{\beta _{1}}}{|C\tilde{e}|^{\alpha _{1}}} \\ 
\leq M_{1}\left\Vert a_{1}\right\Vert _{\infty }((C\mu )^{-\alpha
_{1}}+(C\left\Vert \tilde{e}\right\Vert _{\infty })^{\beta _{1}-\alpha _{1}})
\\ 
\leq C^{\beta _{1}-\alpha _{1}}M_{1}\left\Vert a_{1}\right\Vert _{\infty
}(1+\left\Vert \tilde{e}\right\Vert _{\infty }^{\beta _{1}-\alpha _{1}})%
\text{ \ in }\overline{\Omega },%
\end{array}
\label{6}
\end{equation}%
and%
\begin{equation}
\begin{array}{l}
a_{1}(x)\frac{f_{1}(v)}{(|-C\tilde{e}|+\varepsilon )^{\alpha _{i}}}\geq
-\left\Vert a_{1}\right\Vert _{\infty }\frac{f_{1}(v)}{(|C\tilde{e}%
|+\varepsilon )^{\alpha _{1}}}\geq -M_{1}\left\Vert a_{1}\right\Vert
_{\infty }\frac{1+|v|^{\beta _{1}}}{(|C\tilde{e}|+\varepsilon )^{\alpha _{1}}%
} \\ 
\geq -M_{1}\left\Vert a_{1}\right\Vert _{\infty }\frac{1+|C\tilde{e}|^{\beta
_{1}}}{(|C\tilde{e}|+\varepsilon )^{\alpha _{1}}}\geq -M_{1}\left\Vert
a_{1}\right\Vert _{\infty }\frac{1+|C\tilde{e}|^{\beta _{1}}}{|C\tilde{e}%
|^{\alpha _{1}}} \\ 
\geq -M_{1}\left\Vert a_{1}\right\Vert _{\infty }((C\mu )^{-\alpha
_{1}}+(C\left\Vert \tilde{e}\right\Vert _{\infty })^{\beta _{1}-\alpha _{1}})
\\ 
\geq -C^{\beta _{1}-\alpha _{1}}M_{1}\left\Vert a_{1}\right\Vert _{\infty
}(1+\left\Vert \tilde{e}\right\Vert _{\infty }^{\beta _{1}-\alpha _{1}})%
\text{ \ in }\overline{\Omega },%
\end{array}
\label{6*}
\end{equation}%
provided that $C>0$ is sufficiently large, for all $(u,v)\in \lbrack -C%
\tilde{e},C\tilde{e}]\times \lbrack -C\tilde{e},C\tilde{e}]$ and all $%
\varepsilon \in (0,1)$. Then gathering (\ref{6})-(\ref{6*}) together leads to%
\begin{equation*}
\begin{array}{l}
-\Delta (C\tilde{e})+h_{\lambda ,\phi _{1}}(C\tilde{e})\geq a_{1}(x)\frac{%
f_{1}(v)}{(|C\tilde{e}|+\varepsilon )^{\alpha _{1}}}\text{ \ in }\overline{%
\Omega }%
\end{array}%
\end{equation*}%
and%
\begin{equation*}
\begin{array}{l}
-\Delta (-C\tilde{e})+h_{\lambda ,\phi _{1}}(-C\tilde{e})\leq a_{1}(x)\frac{%
f_{1}(v)}{(|-C\tilde{e}|+\varepsilon )^{\alpha _{1}}}\text{ \ in }\overline{%
\Omega },%
\end{array}%
\end{equation*}%
for all $(u,v)\in \lbrack -C\tilde{e},C\tilde{e}]\times \lbrack -C\tilde{e},C%
\tilde{e}]$ and all $\varepsilon \in (0,1),$ provided that $C>0$ is
sufficiently large. Likewise, a quite similar argument provides%
\begin{equation*}
\begin{array}{l}
-\Delta (C\tilde{e})+h_{\lambda ,\phi _{1}}(C\tilde{e})\geq a_{2}(x)\frac{%
f_{2}(u)}{(|C\tilde{e}|+\varepsilon )^{\alpha _{2}}}\text{ \ in }\overline{%
\Omega }%
\end{array}%
\end{equation*}%
and%
\begin{equation*}
\begin{array}{l}
-\Delta (-C\tilde{e})+h_{\lambda ,\phi _{1}}(-C\tilde{e})\leq a_{2}(x)\frac{%
f_{2}(u)}{(|-C\tilde{e}|+\varepsilon )^{\alpha _{2}}}\text{ \ in }\overline{%
\Omega },%
\end{array}%
\end{equation*}%
for all $(u,v)\in \lbrack -C\tilde{e},C\tilde{e}]\times \lbrack -C\tilde{e},C%
\tilde{e}]$ and all $\varepsilon \in (0,1),$ with $C>0$ is sufficiently
large.

This proves that $(-C\tilde{e},-C\tilde{e})$ and $(C\tilde{e},C\tilde{e})$
are a sub-supersolution pair for $(\mathrm{P}_{\varepsilon })$ for all $%
\varepsilon \in (0,1)$. Consequently, we may apply the general theory of
sub-supersolutions for systems (see, e.g., \cite[Section 5.5]{CLM}) which
ensures the existence of solutions $(u_{\varepsilon },v_{\varepsilon })\in
(H_{0}^{1}(\Omega )\cap L^{\infty }(\Omega ))^{2}$ of $(\mathrm{P}%
_{\varepsilon })$ within $[-C\tilde{e},C\tilde{e}]\times \lbrack -C\tilde{e}%
,C\tilde{e}],$ for all $\varepsilon \in (0,1)$. The proof is now completed.
\end{proof}

\subsection{A sign-changing sub-supersolution pair}

Assume in $\mathrm{H}(a)$ that 
\begin{equation}
\rho _{i}<\frac{1}{2}\underset{\overline{\Omega }}{\max }\phi _{1},
\label{10**}
\end{equation}%
and fix $\gamma _{i}\in (0,1)$ such that%
\begin{equation}
\rho _{i}:=\gamma _{i}^{\frac{-1}{1-\gamma _{i}}},  \label{33}
\end{equation}%
which is possible since $\rho _{i}>1$.

Setting%
\begin{equation}
\overline{u}=\phi _{1}^{\gamma _{1}}-\gamma _{1}\phi _{1},\text{ \ }%
\overline{v}=\phi _{1}^{\gamma _{2}}-\gamma _{2}\phi _{1},  \label{5}
\end{equation}%
observe that 
\begin{equation}
\overline{u}\geq 0\text{ (resp. }\leq 0\text{) if \ }0\leq \phi _{1}\leq
\rho _{1}\text{ (resp. }\phi _{1}\geq \rho _{1}\text{)}  \label{10}
\end{equation}%
and%
\begin{equation}
\overline{v}\geq 0\text{ (resp. }\leq 0\text{) \ if \ }0\leq \phi _{1}\leq
\rho _{2}\text{ (resp. }\phi _{1}\geq \rho _{2}\text{).}  \label{10*}
\end{equation}

\begin{lemma}
\label{L1} Under assumptions $\mathrm{H}(f),$ $\mathrm{H}(a)$ and (\ref{exp}%
), the pair $(\overline{u},\overline{v})$ is a supersolution for problem $(%
\mathrm{P}_{\varepsilon })$, for $\lambda >0$ big enough and for every $%
\varepsilon \in (0,1)$.
\end{lemma}

\begin{proof}
A direct computations shows that%
\begin{equation*}
-\Delta (\phi _{1}^{\gamma _{i}})=\gamma _{i}\lambda _{1}\phi _{1}^{\gamma
_{i}}+\gamma _{i}\left( 1-\gamma _{i}\right) \phi _{1}^{\gamma
_{i}-2}\left\vert \nabla \phi _{1}\right\vert ^{2}\text{ in }\Omega ,\text{
for }i=1,2..
\end{equation*}%
Hence%
\begin{equation}
\begin{array}{l}
-\Delta (\phi _{1}^{\gamma _{i}}-\gamma _{i}\phi _{1})=\gamma _{i}\lambda
_{1}\phi _{1}^{\gamma _{i}}+\gamma _{i}\left( 1-\gamma _{i}\right) \phi
_{1}^{\gamma _{i}-2}\left\vert \nabla \phi _{1}\right\vert ^{2}-\gamma
_{i}\lambda _{1}\phi _{1} \\ 
=\lambda _{1}\gamma _{i}(\phi _{1}^{\gamma _{i}}-\phi _{1})+\gamma
_{i}\left( 1-\gamma _{i}\right) \phi _{1}^{\gamma _{i}-2}\left\vert \nabla
\phi _{1}\right\vert ^{2}\text{ \ in }\Omega ,\text{ for }i=1,2.%
\end{array}
\label{25}
\end{equation}%
We shall prove that $(\overline{u},\overline{v})$\ is a supersolution for
problem $(\mathrm{P}_{\varepsilon })$. To this end, set 
\begin{equation*}
\Omega _{\delta }:=\{x\in \overline{\Omega }:d(x)<\delta \},\text{ with a
constant }\delta >0.
\end{equation*}%
From (\ref{5}), (\ref{10}) and for $\delta >0$ small enough, we have $%
\overline{u}\geq 0,$ $h_{\lambda ,\phi _{1}}(\overline{u})\geq 0$ as well as 
$(\phi _{1}^{\gamma _{i}}-\phi _{1})\geq 0$ in $\Omega _{\delta }$. Thus, by
(\ref{exp}), (\ref{25}), we get%
\begin{equation*}
\begin{array}{l}
(|\overline{u}|+\varepsilon )^{\alpha _{1}}(-\Delta \overline{u}+h_{\lambda
,\phi _{1}}(\overline{u}))\geq (|\phi _{1}^{\gamma _{1}}-\gamma _{1}\phi
_{1}|+\varepsilon )^{\alpha _{1}}(-\Delta \overline{u}) \\ 
\geq \gamma _{1}\left( 1-\gamma _{1}\right) (|\phi _{1}^{\gamma _{1}}-\gamma
_{1}\phi _{1}|+\varepsilon )^{\alpha _{1}}\phi _{1}^{\gamma
_{i}-2}\left\vert \nabla \phi _{1}\right\vert ^{2} \\ 
\geq \gamma _{1}\left( 1-\gamma _{1}\right) (|\phi _{1}^{\gamma
_{1}}(1-\gamma _{1}\phi _{1}^{1-\gamma _{1}})|)^{\alpha _{1}}\phi
_{1}^{\gamma _{i}-2}\left\vert \nabla \phi _{1}\right\vert ^{2} \\ 
\geq \gamma _{1}\left( 1-\gamma _{1}\right) |1-\gamma _{1}\delta ^{1-\gamma
_{1}}|^{\alpha _{1}}\phi _{1}^{\gamma _{1}(1+\alpha _{1})-2}\left\vert
\nabla \phi _{1}\right\vert ^{2}\text{ in }\Omega _{\delta }.%
\end{array}%
\end{equation*}%
Since $0<\gamma _{1}<1$ we have $\gamma _{1}<\frac{2}{1+\alpha _{1}}$ and so 
$\gamma _{1}(1+\alpha _{1})-2<0$ for every $\alpha _{1}\in (0,1)$. Fix $C>0$
such that the conclusion of Theorem \ref{T2} holds true. By (\ref{67}), (\ref%
{19}) and $\mathrm{H}(f)$, we infer that 
\begin{equation*}
\begin{array}{l}
(|\overline{u}|+\varepsilon )^{\alpha _{1}}(-\Delta \overline{u}+h_{\lambda
,\phi _{1}}(\overline{u})) \\ 
\geq \gamma _{1}\left( 1-\gamma _{1}\right) |1-\gamma _{1}\delta ^{1-\gamma
_{1}}|^{\alpha _{1}}(ld(x))^{\gamma _{1}(1+\alpha _{1})-2}\eta ^{2} \\ 
\geq \gamma _{1}\left( 1-\gamma _{1}\right) |1-\gamma _{1}\delta ^{1-\gamma
_{1}}|^{\alpha _{1}}(l\delta )^{\gamma _{1}(1+\alpha _{1})-2}\eta ^{2} \\ 
\geq \left\Vert a_{1}\right\Vert _{\infty }(1+(C\left\Vert \tilde{e}%
\right\Vert _{\infty })^{\beta _{1}})\geq \left\Vert a_{1}\right\Vert
_{\infty }f_{1}(v) \\ 
\geq a_{1}(x)f_{1}(v)\text{ \ in }\Omega _{\delta },%
\end{array}%
\end{equation*}%
for all $v\in \lbrack -C\tilde{e},C\tilde{e}],$ for all $\varepsilon \in
(0,1),$ provided $\delta >0$ is sufficiently small. This shows that%
\begin{equation}
\begin{array}{l}
-\Delta \overline{u}+h_{\lambda ,\phi _{1}}(\overline{u})\geq a_{1}(x)\frac{%
f_{1}(v)}{(|\overline{u}|+\varepsilon )^{\alpha _{1}}}\text{ in }\Omega
_{\delta },%
\end{array}
\label{1}
\end{equation}%
for all $v\in \lbrack -C\tilde{e},C\tilde{e}],$ for all $\varepsilon \in
(0,1)$.

Next, we examine the case when $x\in \Omega _{\rho _{1}}\backslash \overline{%
\Omega }_{\delta }$. From (\ref{24}), (\ref{25}) and (\ref{5}), we have%
\begin{equation}
\begin{array}{l}
-\Delta \overline{u}+h_{\lambda ,\phi _{1}}(\overline{u})\geq \gamma
_{1}\lambda _{1}(\phi _{1}^{\gamma _{1}}-\phi _{1})+\lambda (\phi
_{1}^{\gamma _{1}}+(1-\gamma _{1})\phi _{1}) \\ 
=(\gamma _{1}\lambda _{1}+\lambda )\phi _{1}^{\gamma _{1}}+(\lambda
(1-\gamma _{1})-\gamma _{1}\lambda _{1})\phi _{1})\text{ in }\Omega .%
\end{array}
\label{2}
\end{equation}%
On account of $\mathrm{H}(f)$, (\ref{25}), (\ref{exp}) and recalling that 
\begin{equation*}
\begin{array}{l}
|\phi _{1}^{\gamma _{1}}-\gamma _{1}\phi _{1}|^{\alpha _{1}}>0\text{ \ in }%
\Omega _{\rho _{1}}\backslash \overline{\Omega }_{\delta },%
\end{array}%
\end{equation*}%
we get%
\begin{equation*}
\begin{array}{l}
(|\overline{u}|+\varepsilon )^{\alpha _{1}}(-\Delta \overline{u}+h_{\lambda
,\phi _{1}}(\overline{u})) \\ 
\geq (|\phi _{1}^{\gamma _{1}}-\gamma _{1}\phi _{1}|+\varepsilon )^{\alpha
_{1}}[(\gamma _{1}\lambda _{1}+\lambda )\phi _{1}^{\gamma _{1}}+(\lambda
(1-\gamma _{1})-\gamma _{1}\lambda _{1})\phi _{1})] \\ 
\geq (\gamma _{1}\lambda _{1}+\lambda )|\phi _{1}^{\gamma _{1}}-\gamma
_{1}\phi _{1}|^{\alpha _{1}}\phi _{1}^{\gamma _{1}}\geq \lambda |\phi
_{1}^{\gamma _{1}}-\gamma _{1}\phi _{1}|^{\alpha _{1}}\delta ^{\gamma _{1}}
\\ 
\geq \left\Vert a_{1}\right\Vert _{\infty }(1+(C\left\Vert \tilde{e}%
\right\Vert _{\infty })^{\beta _{1}})\geq a_{1}(x)f_{1}(v)\text{ \ in }%
\Omega _{\rho _{1}}\backslash \overline{\Omega }_{\delta },%
\end{array}%
\end{equation*}

for all $v\in \lbrack -C\tilde{e},C\tilde{e}]$, provided $\lambda >0$ big
enough. Hence, it turns out that%
\begin{equation}
\begin{array}{l}
-\Delta \overline{u}+\lambda \overline{u}\geq a_{1}(x)\frac{f_{1}(v)}{(|%
\overline{u}|+\varepsilon )^{\alpha _{1}}}\text{ in }\Omega _{\rho
_{1}}\backslash \overline{\Omega }_{\delta },%
\end{array}
\label{2*}
\end{equation}%
for all $v\in \lbrack -C\tilde{e},C\tilde{e}],$ for all $\varepsilon \in
(0,1)$.

It remains to prove that the estimate holds true in $\Omega \backslash 
\overline{\Omega }_{\rho _{1}}$. Recall from $\mathrm{H}(a)$ that $a_{1}(x)$
is negative outside $\overline{\Omega }_{\rho _{1}}$. Then, by (\ref{2}), $%
\mathrm{H}(f)$ and for $\lambda >0$ large, it follows that%
\begin{equation*}
\begin{array}{l}
(|\overline{u}|+\varepsilon )^{\alpha _{1}}(-\Delta \overline{u}+h_{\lambda
,\phi _{1}}(\overline{u})) \\ 
\geq |\overline{u}|^{\alpha _{1}}[(\gamma _{1}\lambda _{1}+\lambda )\phi
_{1}^{\gamma _{1}}+(\lambda (1-\gamma _{1})-\gamma _{1}\lambda _{1})\phi
_{1})] \\ 
\geq 0\geq a_{1}(x)f_{1}(v)\text{ \ in }\Omega \backslash \overline{\Omega }%
_{\rho _{1}},%
\end{array}%
\end{equation*}%
for all $v\in \lbrack -C\tilde{e},C\tilde{e}]$ and all $\varepsilon \in
(0,1) $. Thus, it turns out that%
\begin{equation}
\begin{array}{l}
-\Delta \overline{u}+\lambda \overline{u}\geq a_{1}(x)\frac{f_{1}(v)}{(|%
\overline{u}|+\varepsilon )^{\alpha _{1}}}\text{ in }\Omega \backslash 
\overline{\Omega }_{\rho _{1}}.%
\end{array}
\label{0}
\end{equation}%
Gathering together (\ref{1}), (\ref{2*}) and (\ref{0}) we deduce that%
\begin{equation}
\begin{array}{l}
-\Delta \overline{u}+h_{\lambda ,\phi _{1}}(\overline{u})\geq a_{1}(x)\frac{%
f_{1}(v)}{(|\overline{u}|+\varepsilon )^{\alpha _{1}}}\text{ in }\Omega ,%
\end{array}
\label{7}
\end{equation}%
for all $v$ within $[-C\tilde{e},\overline{v}]$, for all $\varepsilon \in
(0,1)$. Similarly, following the same argument as above we obatin%
\begin{equation}
\begin{array}{l}
-\Delta \overline{v}+h_{\lambda ,\phi _{1}}(\overline{v})\geq a_{2}(x)\frac{%
f_{2}(u)}{(|\overline{v}|+\varepsilon )^{\alpha _{2}}}\text{ \ in }\Omega ,%
\end{array}
\label{8}
\end{equation}%
for all $u$ within $[-C\tilde{e},\overline{u}]$, for all $\varepsilon \in
(0,1)$. Consequently, on the basis of (\ref{7}) and (\ref{8}) we conclude
that $(\overline{u},\overline{v})$ is a supersolution of $(\mathrm{P}%
_{\varepsilon })$.
\end{proof}

\begin{remark}
A careful inspection of the proof of Lemma \ref{L1} shows that for a fixed $%
C>0$ in Theorem \ref{T2}, constants $\delta :=\delta (C)$ and $\lambda
:=\lambda (C,\delta )$ can be precisely estimated.
\end{remark}

\begin{lemma}
\label{L2} Under assumptions $\mathrm{H}(f),$ $\mathrm{H}(a)$ and (\ref{exp}%
), problem $(\mathrm{P}_{\varepsilon })$ possesses a subsolution $(%
\underline{u}_{\varepsilon },\underline{v}_{\varepsilon })\in
(H_{0}^{1}(\Omega )\cap L^{\infty }(\Omega ))^{2}$, with $\underline{u}%
_{\varepsilon }\leq \overline{u},$ $\underline{v}_{\varepsilon }\leq 
\overline{v}$, for every $\varepsilon \in (0,1)$ and all $\lambda \geq 0$.
Moreover, there exist functions $\underline{u},\underline{v}\in
H_{0}^{1}(\Omega )$ verifying 
\begin{equation}
\underline{u}<0\text{ in }\Omega \backslash \overline{\Omega }_{\rho _{1}},%
\text{ \ \ }\underline{u}\geq 0\text{ \ in }\Omega _{\rho _{1}}  \label{14*}
\end{equation}%
and%
\begin{equation}
\underline{v}<0\text{ in }\Omega \backslash \overline{\Omega }_{\rho _{2}},%
\text{ \ \ }\underline{v}\geq 0\text{ \ in }\Omega _{\rho _{2}},
\label{14**}
\end{equation}
\ such that%
\begin{equation}
\underline{u}_{\varepsilon }\rightarrow \underline{u}\ \ \text{and }\ 
\underline{v}_{\varepsilon }\rightarrow \underline{v}\text{ \ in }%
H_{0}^{1}(\Omega )\text{ as }\varepsilon \rightarrow 0.  \label{14}
\end{equation}
\end{lemma}

\begin{proof}
For any $s\in 
\mathbb{R}
$, denote by $s_{+}:=\max \{s,0\}$ and $s_{-}:=\max \{-s,0\}$.

Define the truncation 
\begin{equation}
\chi _{\phi _{1}}(s)=\frac{1}{\left\Vert \phi _{1}\right\Vert _{\infty }}%
\left\{ 
\begin{array}{ll}
\phi _{1} & \text{if }s\geq 2\phi _{1} \\ 
s-\phi _{1} & \text{if }\phi _{1}\leq s\leq 2\phi _{1} \\ 
0 & \text{if }s\leq \phi _{1}%
\end{array}%
\right.  \label{11}
\end{equation}%
and consider the problem 
\begin{equation*}
(\mathrm{\tilde{P}}_{\varepsilon })\qquad \left\{ 
\begin{array}{ll}
-\Delta u+h_{\lambda ,\phi _{1}}(u)=\mathcal{F}_{1,\varepsilon }(x,u,v) & 
\text{in }\Omega \\ 
-\Delta v+h_{\lambda ,\phi _{1}}(v)=\mathcal{F}_{2,\varepsilon }(x,u,v) & 
\text{in }\Omega \\ 
u,v=0 & \text{on }\partial \Omega ,%
\end{array}%
\right.
\end{equation*}%
for $\varepsilon \in (0,1)$, where%
\begin{equation}
\mathcal{F}_{1,\varepsilon }(x,u,v)=\left\{ 
\begin{array}{ll}
a_{1,+}(x)\chi _{\phi _{1}}(u_{+})\frac{f_{1}(v)}{(|\overline{u}|+1)^{\alpha
_{1}}} & \text{in }\Omega _{\rho _{1}} \\ 
-a_{1,-}(x)\frac{1+|\overline{v}|^{\beta _{1}}}{(|u|+\varepsilon )^{\alpha
_{1}}} & \text{in }\Omega \backslash \overline{\Omega }_{\rho _{1}}%
\end{array}%
\right.  \label{16}
\end{equation}%
and%
\begin{equation}
\mathcal{F}_{2,\varepsilon }(x,u,v)=\left\{ 
\begin{array}{ll}
a_{2,+}(x)\chi _{\phi _{1}}(v_{+})\frac{f_{2}(u)}{(|\overline{v}|+1)^{\alpha
_{2}}} & \text{in }\Omega _{\rho _{2}} \\ 
-a_{2,-}(x)\frac{1+|\overline{u}|^{\beta _{1}}}{(|v|+\varepsilon )^{\alpha
_{2}}} & \text{in }\Omega \backslash \overline{\Omega }_{\rho _{2}}.%
\end{array}%
\right.  \label{16*}
\end{equation}%
It is a simple matter to see that 
\begin{equation}
\begin{array}{l}
\mathcal{F}_{1,\varepsilon }(x,u,v)\leq a_{1}(x)\frac{f_{1}(v)}{%
(|u|+\varepsilon )^{\alpha _{1}}}\text{ \ and \ }\mathcal{F}_{2,\varepsilon
}(x,u,v)\leq a_{2}(x)\frac{f_{2}(v)}{(|v|+\varepsilon )^{\alpha _{2}}},%
\end{array}
\label{31}
\end{equation}%
for all $(u,v)\in \lbrack -C\tilde{e},C\tilde{e}]\times \lbrack -C\tilde{e},C%
\tilde{e}]$ and all $\varepsilon \in (0,1)$. Then, any solution of $(\mathrm{%
\tilde{P}}_{\varepsilon })$ within $[-C\tilde{e},C\tilde{e}]\times \lbrack -C%
\tilde{e},C\tilde{e}]$ is a subsolution of $(\mathrm{P}_{\varepsilon })$.

We claim that $(-C\tilde{e},-C\tilde{e})$ is a subsolution of $(\mathrm{%
\tilde{P}}_{\varepsilon })$. Indeed, from (\ref{62}), (\ref{18}), (\ref{3}), 
$\mathrm{H}(a)$ and $\mathrm{H}(f)$, we have%
\begin{equation*}
\begin{array}{l}
-\Delta (-C\tilde{e})+h_{\lambda ,\phi _{1}}(-C\tilde{e})\leq 0\leq \tilde{f}%
_{\varepsilon }(x,-C\tilde{e},v)\text{ \ in }\Omega _{\rho _{1}}%
\end{array}%
\end{equation*}%
and%
\begin{equation*}
\begin{array}{l}
-\Delta (-C\tilde{e})+h_{\lambda ,\phi _{1}}(-C\tilde{e})\leq -C(1+\lambda 
\frac{\mu }{c})+\lambda \left\Vert \phi _{1}\right\Vert _{\infty } \\ 
\leq -C^{\alpha _{1}}\left\Vert a_{1,-}\right\Vert _{\infty }\mu ^{\alpha
_{1}}\left\Vert \overline{v}\right\Vert _{\infty }^{\beta _{1}}\leq
-C^{-\alpha _{1}}\left\Vert a_{1,-}\right\Vert _{\infty }\frac{1+\left\Vert 
\overline{v}\right\Vert _{\infty }^{\beta _{1}}}{(c^{-1}\mu )^{\alpha _{1}}}
\\ 
\leq -a_{1,-}(x)\frac{1+|\overline{v}|^{\beta _{1}}}{(|-C\tilde{e}%
|+\varepsilon )^{\alpha _{1}}}=\mathcal{F}_{1,\varepsilon }(x,-C\tilde{e},v)%
\text{ \ in }\Omega \backslash \overline{\Omega }_{\rho _{1}},%
\end{array}%
\end{equation*}%
provided that $C>0$ is sufficiently large, for all $v\in \lbrack -C\tilde{e},%
\overline{v}]$, and all $\varepsilon \in (0,1)$. Similarly, one derives that%
\begin{equation*}
-\Delta (-C\tilde{e})+h_{\lambda ,\phi _{1}}(-C\tilde{e})\leq \mathcal{F}%
_{2,\varepsilon }(x,u,-C\tilde{e})\text{ \ in }\overline{\Omega },
\end{equation*}%
provided that $C>0$ is sufficiently large, for all $u\in \lbrack -C\tilde{e},%
\overline{u}]$ and all $\varepsilon \in (0,1)$. This proves the claim.

On the basis of (\ref{31}) and Lemma \ref{L1}, $(\overline{u},\overline{v})$
in (\ref{5}) is a supersolution of problem $(\mathrm{\tilde{P}}_{\varepsilon
})$. Consequently, owing to \cite[section 5.5]{CLM}, problem $(\mathrm{%
\tilde{P}}_{\varepsilon })$ admits a solution $(\underline{u}_{\varepsilon },%
\underline{v}_{\varepsilon })$ within $[-Ce,\overline{u}]\times \lbrack -Ce,%
\overline{v}]$ and $\underline{u}_{\varepsilon },\underline{v}_{\varepsilon
} $ are both nontrivial because $\overline{u}$ and $\overline{v}$ are both
of sign-changing in $\Omega $. Moreover, according to (\ref{31}), $(%
\underline{u}_{\varepsilon },\underline{v}_{\varepsilon })$ is a subsolution
of $(\mathrm{P}_{\varepsilon })$ for all $\varepsilon \in (0,1)$.

Now we prove (\ref{14}). Set $\varepsilon =\frac{1}{n}$ with any positive
integer $n>1$. From above there exist $\underline{u}_{n}:=\underline{u}_{%
\frac{1}{n}}$ and $\underline{v}_{n}:=\underline{v}_{\frac{1}{n}}$ such that%
\begin{equation}
\left\{ 
\begin{array}{c}
\langle -\Delta \underline{u}_{n}+h_{\lambda ,\phi _{1}}(\underline{u}%
_{n}),\varphi \rangle =\int_{\Omega }\mathcal{F}_{1,n}(x,\underline{u}_{n},%
\underline{v}_{n})\varphi \ dx \\ 
\langle -\Delta \underline{v}_{n}+h_{\lambda ,\phi _{1}}(\underline{v}%
_{n}),\psi \rangle =\int_{\Omega }\mathcal{F}_{2,n}(x,\underline{u}_{n},%
\underline{v}_{n})\psi \ dx%
\end{array}%
\right.  \label{22*}
\end{equation}%
for all $\varphi ,\psi \in H_{0}^{1}(\Omega )$ with%
\begin{equation}
-C\tilde{e}\leq \underline{u}_{n}\leq \overline{u}\leq C\tilde{e}\text{, \ }%
-C\tilde{e}\leq \underline{v}_{n}\leq \overline{v}\leq C\tilde{e}\ \text{ in 
}\Omega .  \label{23}
\end{equation}%
Acting with $\varphi =\underline{u}_{n}$ in (\ref{22*}), by (\ref{16}), $%
\mathrm{H}(a)$ and (\ref{exp}), bearing in mind (\ref{23}), it follows that 
\begin{equation*}
\begin{array}{l}
\int_{\Omega }(|\nabla \underline{u}_{n}|^{2}+h_{\lambda ,\phi _{1}}(%
\underline{u}_{n}))\ dx=\int_{\Omega }(|\nabla \underline{u}%
_{n}|^{2}+\lambda (|\underline{u}_{n}|^{2}+\phi _{1}\underline{u}_{n}))\ dx
\\ 
=\int_{\Omega _{\rho _{1}}}a_{1,+}(x)\chi _{\phi _{1}}(\underline{u}_{n,+})%
\frac{f_{1}(\underline{v}_{n})}{(|\overline{u}|+1)^{\alpha _{1}}}\underline{u%
}_{n}\ dx-\int_{\Omega \backslash \overline{\Omega }_{\rho _{1}}}a_{1,-}(x)%
\frac{1+|\overline{v}|^{\beta _{1}}}{(|\underline{u}_{n}|+\varepsilon
)^{\alpha _{1}}}\underline{u}_{n}\ dx \\ 
\leq \int_{\Omega _{\rho _{1}}}\left\Vert a_{1}\right\Vert _{\infty }(1+(C%
\tilde{e})^{\beta _{1}})\underline{u}_{n}\ dx-\int_{\Omega \backslash 
\overline{\Omega }_{\rho _{1}}}a_{1,-}(x)\frac{1}{(C\tilde{e}+1)^{\alpha
_{1}}}\underline{u}_{n}\ dx \\ 
\leq |\Omega |\left\Vert a_{1}\right\Vert _{\infty }(1+C\left\Vert \tilde{e}%
\right\Vert _{\infty }^{\beta _{1}})C\left\Vert \tilde{e}\right\Vert
_{\infty }<\infty .%
\end{array}%
\end{equation*}%
This proves that $\{\underline{u}_{n}\}$ is bounded in $H_{0}^{1}\left(
\Omega \right) $. Similarly, we derive that $\{\underline{v}_{n}\}$ is
bounded in $H_{0}^{1}\left( \Omega \right) $. We are thus allowed to extract
subsequences (still denoted by $\{\underline{u}_{n}\}$ and $\{\underline{v}%
_{n}\}$ ) such that 
\begin{equation*}
\begin{array}{c}
\underline{u}_{n}\rightharpoonup \underline{u}\text{ \ and }\underline{v}%
_{n}\rightharpoonup \underline{v}\text{ in }H_{0}^{1}(\Omega )%
\end{array}%
\end{equation*}%
and 
\begin{equation*}
-C\tilde{e}\leq \underline{u}\leq \overline{u}\leq C\tilde{e}\text{, \ }-C%
\tilde{e}\leq \underline{v}\leq \overline{v}\leq C\tilde{e}\ \text{ in }%
\Omega .
\end{equation*}%
Inserting $(\varphi ,\psi )=(\underline{u}_{n}-\underline{u},\underline{v}%
_{n}-\underline{v})$ in (\ref{22*}) yields \ 
\begin{equation*}
\left\{ 
\begin{array}{c}
\langle -\Delta \underline{u}_{n}+h_{\lambda ,\phi _{1}}(\underline{u}_{n}),%
\underline{u}_{n}-\underline{u}\rangle =\int_{\Omega }\mathcal{F}_{1,n}(x,%
\underline{u}_{n},\underline{v}_{n})(\underline{u}_{n}-\underline{u})\ dx \\ 
\langle -\Delta \underline{v}_{n}+h_{\lambda ,\phi _{1}}(\underline{v}_{n}),%
\underline{v}_{n}-\underline{v}\rangle =\int_{\Omega }\mathcal{F}_{2,n}(x,%
\underline{u}_{n},\underline{v}_{n})(\underline{v}_{n}-\underline{v})\ dx.%
\end{array}%
\right.
\end{equation*}%
By (\ref{13}), (\ref{23}), $\mathrm{H}(a)$ and (\ref{exp}), we have%
\begin{equation*}
\begin{array}{l}
|a_{1,+}(x)\chi _{\phi _{1}}(\underline{u}_{n,+})\frac{f_{1}(\underline{v}%
_{n})}{(|\overline{u}|+1)^{\alpha _{1}}}(\underline{u}_{n}-\underline{u}%
)|\leq 2\left\Vert a_{1}\right\Vert _{\infty }(1+C\left\Vert \tilde{e}%
\right\Vert _{\infty }^{\beta _{1}})C\left\Vert \tilde{e}\right\Vert
_{\infty }%
\end{array}%
\end{equation*}%
as well as%
\begin{equation*}
\begin{array}{l}
|a_{1,-}(x)\frac{1+|\overline{v}|^{\beta _{1}}}{(|\underline{u}%
_{n}|+\varepsilon )^{\alpha _{1}}}(\underline{u}_{n}-\underline{u})|\leq
\left\Vert a_{1}\right\Vert _{\infty }|\underline{u}_{n}|^{1-\alpha _{1}}(1+|%
\overline{v}|^{\beta _{1}}) \\ 
\leq \left\Vert \Lambda _{1}\right\Vert _{\infty }(C\left\Vert \tilde{e}%
\right\Vert _{\infty })^{1+\alpha _{1}}(1+(C\left\Vert \tilde{e}\right\Vert
_{\infty })^{\beta _{1}}).%
\end{array}%
\end{equation*}%
Then, applying Fatou's Lemma, it follows that%
\begin{equation*}
\begin{array}{l}
\underset{n\rightarrow \infty }{\limsup }\int_{\Omega }\mathcal{F}_{1,n}(x,%
\underline{u}_{n},\underline{v}_{n})(\underline{u}_{n}-\underline{u})\ dx \\ 
\leq \int_{\Omega }\underset{n\rightarrow \infty }{\lim }\sup (\mathcal{F}%
_{1,n}(x,\underline{u}_{n},\underline{v}_{n})(\underline{u}_{n}-\underline{u}%
))\ dx\rightarrow 0\text{ as }n\rightarrow +\infty ,%
\end{array}%
\end{equation*}%
showing that $\underset{n\rightarrow \infty }{\lim }\left\langle -\Delta 
\underline{u}_{n},\underline{u}_{n}-\underline{u}\right\rangle \leq 0.$ Then
the $S_{+}$-property of $-\Delta $ on $H_{0}^{1}(\Omega )$ guarantees that $%
\underline{u}_{n}\rightarrow \underline{u}$ in $H_{0}^{1}(\Omega )$.
Similarly, we prove that $\underline{v}_{n}\rightarrow \underline{v}$ in $%
H_{0}^{1}(\Omega )$. Hence, we may pass to the limit in (\ref{22*}) to
conclude that $(\underline{u},\underline{v})$ is a solution of problem $(%
\mathrm{\tilde{P}}_{0})$ (That is $(\mathrm{\tilde{P}}_{\varepsilon })$ with 
$\varepsilon =0$).

It remains to prove (\ref{14*}) and (\ref{14**}). Since $\underline{u}\leq 
\overline{u}$ and $\overline{u}$ is negative in $\Omega \backslash \overline{%
\Omega }_{\rho _{1}}$, it turns out that $\underline{u}$ is negative in $%
\Omega \backslash \overline{\Omega }_{\rho _{1}}$. Similarly, one derives
that $\underline{v}$ is negative in $\Omega \backslash \overline{\Omega }%
_{\rho _{2}}$. Acting in $(\mathrm{\tilde{P}}_{0})$ with $(1_{\Omega _{\rho
_{1}}}\underline{u}_{-},1_{\Omega _{\rho _{1}}}\underline{v}_{-}),$ from the
definition of the truncation $\chi _{\phi _{1}}$ in (\ref{11}), it follows
that%
\begin{equation*}
\begin{array}{l}
\int_{\Omega _{\rho _{1}}}(|\nabla \underline{u}_{-}|^{2}+\lambda (|%
\underline{u}_{-}|^{2}+\phi _{1}\underline{u}_{-}))\text{ }dx \\ 
=\int_{\Omega _{\rho _{1}}}a_{1,+}(x)\chi _{\phi _{1}}(u_{+})\frac{f_{1}(v)}{%
(|\overline{u}|+1)^{\alpha _{1}}}\underline{u}_{-}\text{ }dx=0%
\end{array}%
\end{equation*}%
and%
\begin{equation*}
\begin{array}{l}
\int_{\Omega _{\rho _{1}}}(|\nabla \underline{v}_{-}|^{2}+\lambda (|%
\underline{v}_{-}|^{2}+\phi _{1}\underline{v}_{-}))\text{ }dx \\ 
=\int_{\Omega _{\rho _{2}}}\Lambda _{2,+}(x)\chi _{\phi _{1}}(v)|u|^{\alpha
_{2}}(v_{+}+1)^{\beta _{2}}\underline{v}_{-}\text{ }dx=0.%
\end{array}%
\end{equation*}%
Consequently, one can conclude that $\underline{u}\geq 0$ in $\Omega _{\rho
_{1}}$ and $\underline{v}\geq 0$ in $\Omega _{\rho _{2}}$. This completes
the proof.
\end{proof}

\section{Existence of solutions}

\label{S3}

Our main result is formulated as follows.

\begin{theorem}
\label{T1} Under assumptions $\mathrm{H}(f),$ $\mathrm{H}(a)$ and (\ref{exp}%
), problem $(\mathrm{P})$ possesses a nodal solution $(\breve{u},\breve{v}%
)\in (H_{0}^{1}(\Omega )\cap L^{\infty }(\Omega ))^{2}$ for $\lambda >0$
large.
\end{theorem}

\begin{proof}
Set $\varepsilon =\frac{1}{n}$ with any positive integer $n>1$. According to
Lemmas \ref{L1} and \ref{L2} and thanks to \cite[Section 5.5]{CLM}, there
exists $(u_{n},v_{n}):=(u_{\frac{1}{n}},v_{\frac{1}{n}})\in
(H_{0}^{1}(\Omega )\cap L^{\infty }(\Omega ))^{2}$ such that%
\begin{equation}
\left\{ 
\begin{array}{c}
\langle -\Delta u_{n}+h_{\lambda ,\phi _{1}}(u_{n}),\varphi \rangle
=\int_{\Omega }a_{1}(x)\frac{f_{1}(v_{n})}{(|u_{n}|+\frac{1}{n})^{\alpha
_{1}}}\varphi \ dx \\ 
\langle -\Delta v_{n}+h_{\lambda ,\phi _{1}}(v_{n}),\psi \rangle
=\int_{\Omega }a_{2}(x)\frac{f_{2}(u_{n})}{(|v_{n}|+\frac{1}{n})^{\alpha
_{2}}}\psi \ dx%
\end{array}%
\right.  \label{22}
\end{equation}%
for all $\varphi ,\psi \in H_{0}^{1}(\Omega )$ and%
\begin{equation}
\underline{u}_{n}\leq u_{n}\leq \overline{u}\text{, }\underline{v}_{n}\leq
v_{n}\leq \overline{v}\text{, }\forall n\text{.}  \label{15}
\end{equation}%
In addition, due to (\ref{14*}) and (\ref{14**}) one has%
\begin{equation}
u_{n}<0\text{ in }\Omega \backslash \overline{\Omega }_{\rho _{1}},\text{ \
\ }u_{n}\geq 0\text{ \ in }\Omega _{\rho _{1}}
\end{equation}%
and%
\begin{equation}
v_{n}<0\text{ in }\Omega \backslash \overline{\Omega }_{\rho _{2}},\text{ \
\ }v_{n}\geq 0\text{ \ in }\Omega _{\rho _{2}},\text{ }\forall n.
\end{equation}%
We claim that the sets where the sequences $\{u_{n}\}$ and $\{v_{n}\}$
vanish are negligible. To this end, define the measurable set 
\begin{equation*}
\Gamma _{n}:=\{x\in \Omega _{\rho _{1}}:u_{n}(x)=0\},\text{ }\forall n.
\end{equation*}%
By a classical result of measure theory (see \cite[Theorem 3.28]{WZ}), $%
\Gamma _{n}$ is a $\mathcal{G}_{\gamma }$-set for less than a zero measure
set. Thus, one can write 
\begin{equation}
\Gamma _{n}=\mathcal{A}_{n}-\mathcal{B}_{n},  \label{77}
\end{equation}%
where $\mathcal{A}_{n}-\mathcal{B}_{n}$ is the relative complement of $%
\mathcal{A}_{n}$ in $\mathcal{B}_{n},$ $|\mathcal{B}_{n}|=0$ and $\mathcal{A}%
_{n}$ is a $\mathcal{G}_{\gamma }$-set, that is 
\begin{equation*}
\mathcal{A}_{n}=\underset{k}{\cap }\mathcal{G}_{n,k},\text{ \ }\mathcal{G}%
_{n,k}\text{ is open.}
\end{equation*}%
Without loss of generality, we may assume that $\mathcal{B}_{n}$ is not
dense in $\Omega _{\rho _{1}}$ because otherwise, by density and according
to \cite[Proposition 2.4]{BC}, the term $a_{1}(x)\frac{f_{1}(v_{n})}{%
(|u_{n}|+\frac{1}{n})^{\alpha _{1}}}$ should vanish in $\Omega _{\rho _{1}}$
which is absurd in view of $\mathrm{H}(f)$ and $\mathrm{H}(a)$.

Let $\tilde{\varphi}\in C_{0}^{\infty }(\mathcal{A}_{n})$ such that 
\begin{equation*}
\tilde{\varphi}>0\text{ in }\mathcal{A}_{n}\text{ \ and \ }\tilde{\varphi}=0%
\text{ in }\Omega _{\rho _{1}}\backslash \mathcal{A}_{n}.
\end{equation*}%
Testing with $\tilde{\varphi}$ in (\ref{22}) we have%
\begin{equation*}
\int_{\mathcal{A}_{n}}(\nabla u_{n}\nabla \tilde{\varphi}+h_{\lambda ,\phi
_{1}}(u_{n})\tilde{\varphi})\text{ }\mathrm{d}x=\int_{\mathcal{A}%
_{n}}a_{1}(x)\frac{f_{1}(v_{n})}{(|u_{n}|+\frac{1}{n})^{\alpha _{1}}}\tilde{%
\varphi}\text{ }\mathrm{d}x
\end{equation*}%
In view of (\ref{77}) we get%
\begin{equation}
\int_{\Gamma _{n}}(\nabla (u_{n}\nabla \tilde{\varphi}+h_{\lambda ,\phi
_{1}}(u_{n})\tilde{\varphi})\text{ }\mathrm{d}x=\int_{\Gamma _{n}}a_{1}(x)%
\frac{f_{1}(v_{n})}{(|u_{n}|+\frac{1}{n})^{\alpha _{1}}}\tilde{\varphi}\text{
}\mathrm{d}x.  \label{55}
\end{equation}%
Since $u_{n}\in W_{loc}^{1,1}(\Omega )$ it follows, by \cite[Lemma 7.7]{GT},
that $\nabla u_{n}=0$ on $\Gamma _{n}$. Replacing in (\ref{55}) leads%
\begin{equation*}
\int_{\Gamma _{n}}\lambda \phi _{1}\tilde{\varphi}\text{ }\mathrm{d}%
x=n^{\alpha _{1}}\int_{\Gamma _{n}}a_{1}(x)f_{1}(v_{n})\tilde{\varphi}\text{ 
}dx\text{, }\forall n\geq 1.
\end{equation*}%
which, by $\mathrm{H}(f),$ $\mathrm{H}(a),$ forces that $|\Gamma _{n}|=0$.
The same conclusion can be drawn for the set $\hat{\Gamma}_{n}:=\{x\in
\Omega _{\rho _{2}}:v_{n}(x)=0\},$ for all $n\geq 1$. This shows that $u_{n}$
and $v_{n}$ cannot be identically zero in $\Omega _{\rho _{1}}$ and $\Omega
_{\rho _{2}}$, respectively.

Now, taking $\varphi =u_{n}$ in (\ref{22}), by $\mathrm{H}(f),$ $\mathrm{H}%
(a)$ and since%
\begin{equation}
|u_{n}|,|v_{n}|\leq C\tilde{e}\ \text{ in }\Omega ,  \label{13}
\end{equation}%
we get 
\begin{equation}
\begin{array}{l}
\int_{\Omega }(|\nabla u_{n}|^{2}+h_{\lambda ,\phi _{1}}(u_{n})u_{n})\
dx\leq \int_{\Omega }\left\Vert a_{1}\right\Vert _{\infty }\frac{%
1+|v_{n}|^{\beta _{1}}}{(|u_{n}|+\frac{1}{n})^{\alpha _{1}}}u_{n}\ dx \\ 
\leq \left\Vert a_{1}\right\Vert _{\infty }\int_{\Omega }|u_{n}|^{1-\alpha
_{1}}(1+|v_{n}|^{\beta _{1}})\ dx \\ 
\leq \left\Vert a_{1}\right\Vert _{\infty }\int_{\Omega }(C\tilde{e}%
)^{1-\alpha _{1}}(1+(C\tilde{e})^{\beta _{1}})\ dx<\infty ,%
\end{array}
\label{123}
\end{equation}%
showing that $\{u_{n}\}$ is bounded in $H_{0}^{1}\left( \Omega \right) $.
Similarly, we prove that $\{v_{n}\}$ is bounded in $H_{0}^{1}\left( \Omega
\right) $. We are thus allowed to extract subsequences (still denoted by $%
\{u_{n}\}$ and $\{v_{n}\}$) such that 
\begin{equation}
\begin{array}{c}
u_{n}\rightharpoonup \breve{u}\text{ \ and }v_{n}\rightharpoonup \breve{v}%
\text{ in }H_{0}^{1}(\Omega )%
\end{array}
\label{30}
\end{equation}%
Moreover, by (\ref{14}) and (\ref{15}) one has 
\begin{equation}
\underline{u}\leq \breve{u}\leq \overline{u}\text{, \ \ }\underline{v}\leq 
\breve{v}\leq \overline{v}\ \text{ in }\Omega .  \label{12}
\end{equation}%
Then, analysis similar to that of the proof of (\ref{14}) in Lemma \ref{L2}
leads to 
\begin{equation*}
\begin{array}{c}
u_{n}\rightarrow \breve{u}\text{ \ and }v_{n}\rightarrow \breve{v}\text{ in }%
H_{0}^{1}(\Omega ).%
\end{array}%
\end{equation*}%
Setting%
\begin{equation*}
\Gamma :=\{x\in \Omega _{\rho _{1}}:\breve{u}(x)=0\}\text{ and }\breve{\Gamma%
}:=\{x\in \Omega _{\rho _{2}}:\breve{v}(x)=0\},
\end{equation*}%
by (\ref{13}) and the dominated convergence theorem, it follows that 
\begin{equation}
|\Gamma |=\lim_{n\rightarrow +\infty }|\Gamma _{n}|=0\text{ \ and \ }|\breve{%
\Gamma}|=\lim_{n\rightarrow +\infty }|\breve{\Gamma}_{n}|=0,  \label{51}
\end{equation}%
showing that $\breve{u}$ and $\breve{v}$ cannot be identically zero in $%
\Omega _{\rho _{1}}$ and $\Omega _{\rho _{2}}$, respectively.

Finally, we may pass to the limit in (\ref{22}) to conclude that $(\breve{u},%
\breve{v})$ is a solution of problem $(\mathrm{P})$ within $[\underline{u},%
\overline{u}]\times \lbrack \underline{v},\overline{v}]$. On account of (\ref%
{10**})-(\ref{10*}) and (\ref{51}) together with Lemma \ref{L2}, we infer
that $\breve{u},\breve{v}$ are both of sign-changing and satisfying%
\begin{equation*}
\breve{u}\leq \overline{u}<0\text{ in }\Omega \backslash \overline{\Omega }%
_{\rho _{1}}\text{ and \ }\breve{u}\geq \underline{u}\geq 0,\text{ }\breve{u}%
\neq 0\text{ in }\Omega _{\rho _{1}},
\end{equation*}%
\begin{equation*}
\breve{v}\leq \overline{v}<0\text{ in }\Omega \backslash \overline{\Omega }%
_{\rho _{2}}\text{ and \ }\breve{v}\geq \underline{v}\geq 0,\text{ }\breve{v}%
\neq 0\text{ in }\Omega _{\rho _{2}}.
\end{equation*}%
This completes the proof.
\end{proof}

\begin{acknowledgement}
The author was supported by the Directorate-General of Scientific Research
and Technological Development (DGRSDT).
\end{acknowledgement}

\end{document}